\documentclass[onefignum,onetabnum]{siamart250211}

\usepackage{lipsum}
\usepackage{amsfonts}
\usepackage{graphicx}
\usepackage{epstopdf}
\usepackage{algorithmic}
\ifpdf
  \DeclareGraphicsExtensions{.eps,.pdf,.png,.jpg}
\else
  \DeclareGraphicsExtensions{.eps}
\fi

\newsiamremark{remark}{Remark}
\newsiamremark{hypothesis}{Hypothesis}
\crefname{hypothesis}{Hypothesis}{Hypotheses}
\newsiamthm{claim}{Claim}

\usepackage[most]{tcolorbox}
\usepackage{xcolor}
\usepackage{csquotes}

\definecolor{summarybg}{HTML}{F9F9F9}
\definecolor{summaryborder}{HTML}{B0B0B0}

\newtcolorbox{summarybox}[1][Summary]{
  colback=summarybg,
  colframe=summaryborder,
  coltitle=black,
  fonttitle=\bfseries,
  title=#1,
  boxrule=0.8pt,
  arc=4pt,
  left=6pt,
  right=6pt,
  top=4pt,
  bottom=4pt,
  enhanced
}

\headers{Differentiating through Stochastic Differential Equations}{R. Leburu, L. Nurbekyan, L. Ruthotto}

\title{Differentiating through Stochastic Differential Equations: A Primer\thanks{Submitted to the editors January 14, 2026.
\funding{This work was supported in part by NSF award DMS 2038118.}}}

\author{
  Rishi Leburu\thanks{Department of Mathematics, Emory University, Atlanta, GA 
  (\email{rishileburu@gmail.com}, \email{lnurbek@emory.edu}).}
  \and 
  Levon Nurbekyan\footnotemark[2]
  \and 
  Lars Ruthotto\thanks{Departments of Mathematics and Computer Science, Emory University, Atlanta, GA, USA (\email{lruthotto@emory.edu}, \url{https://www.math.emory.edu/\string~lruthot/})}
}
\usepackage{amsopn}

\usepackage{booktabs}
\usepackage{placeins}

\ifpdf
\hypersetup{
  pdftitle={Differentiating through Stochastic Differential Equations},
  pdfauthor={R. Leburu, L. Nurbekyan, and L. Ruthotto}
}
\fi

\begin{document}

\maketitle

\begin{abstract}
    Dynamical systems are essential to model various phenomena in  physics, finance, economics, and are also of current interest in machine learning. A central modeling task is investigating parameter sensitivity, whether tuning atmospheric coefficients, computing financial Greeks, or optimizing neural networks. These sensitivities are mathematically expressed as derivatives of an objective function with respect to parameters of interest and are rarely available analytically, necessitating numerical methods for approximating them. While the literature for differentiation of deterministic systems is well-covered, the treatment of stochastic systems, such as stochastic differential equations (SDEs), in most curricula is less comprehensive than what the subtleties arising from the interplay of noise and discretization warrant.
    
    This paper provides a primer on numerical differentiation of SDEs organized as a two-tale narrative. Tale 1 demonstrates that differentiating through discretized SDEs, known as the discretize-optimize approach, is reliable for both It\^{o} and Stratonovich calculus. Tale 2 examines the optimize-discretize approach, investigating the continuous limit of adjoint equations from Tale 1 corresponding to the desired gradients. Our aim is to equip readers with a clear guide on the numerical differentiation of SDEs: computing gradients correctly in both It\^{o} and Stratonovich settings, understanding when discretize-optimize and optimize-discretize agree or diverge, and developing intuition for reasoning about stochastic differentiation beyond the cases explicitly covered.
\end{abstract}



\section{Introduction} 

Dynamical systems are a natural modeling framework for various applications in physics, economics, finance, and have recently also been used in machine learning. Familiar examples are the Lorenz-63 model for atmospheric convection~\cite{lorenz1963}, the Black--Scholes model~\cite{black1973,merton1973} for the evolution of a stock price in financial markets, and neural ordinary differential equations~\cite{chen2018} (neural ODEs) for modeling complex nonlinear input-output transformations in machine learning.

A key part of modeling with dynamical systems is investigating the sensitivity of these systems with respect to parameters that define their dynamics or initial conditions. For Lorenz-63, for example, one is interested in tuning the system's coefficients to match meteorological data. For Black--Scholes models, one computes sensitivities to market parameters such as volatility, interest rate, and maturity, as well as the sensitivity of the future stock price with respect to the current price; quantities collectively known as Greeks~\cite[Ch.~19]{hull2018}. Finally, for neural ODEs, one seeks the weights of the networks that define the dynamics that yield the best predictive models.

Mathematically, these sensitivities can be expressed as derivatives of an objective function with respect to the parameters of interest. Their computation involves a secondary dynamical system whose solution encodes how the primary system responds to perturbations. As the solution is rarely available analytically, they are usually approximated with numerical methods. Such methods take one of two natural routes: differentiating a discretization of the primary system, or first deriving the secondary continuous-time system and then discretizing its solution. Thus, the main goal of this paper is to introduce the reader to the basic mathematical and numerical principles underpinning the computation of these sensitivities.

In this paper, our aim is to cover the differentiation of stochastic differential equations for the following primary reasons.

\begin{enumerate}
    \item The literature for deterministic systems is very rich, and there are numerous excellent introductory resources, including textbooks~\cite[Chap.~5]{kirk2004}. Indeed, computing sensitivities of deterministic ODEs with respect to input parameters is part of the standard material in optimal control (OC) theory.
    \item The literature on stochastic systems is more limited. Often, standard textbooks on stochastic differential equations (SDEs)~\cite{ higham2021, kloeden1992, oksendal2003} or stochastic optimal control~\cite{ fleming2006, pham2009, yong1999} do not include discussions of sensitivities akin to those in ODEs.
    \item Stochastic systems are more delicate to work with, given various modeling and approximation considerations, such as noise (stochasticity) type and numerical discretization. This delicacy leads to interesting mathematical subtleties that are often overlooked and worth discussing explicitly.
    \item Differentiating through SDEs arises naturally in several important applications, for example: computing option sensitivities (Greeks) in quantitative finance~\cite[Chap.~19]{hull2018}, optimizing control policies for systems governed by stochastic dynamics~\cite[Chap.~6]{pham2009}, and training neural SDE models in machine learning~\cite{kidger2022neural}. These tasks all require computing gradients of objectives defined through stochastic dynamics, motivating the techniques discussed in this paper.
\end{enumerate}

This tutorial serves as an accessible introduction to SDE differentiation, geared toward graduate or advanced undergraduates, faculty, and practitioners who are new to it. As a classroom-style tutorial, this is not a survey covering state-of-the-art 
methods. We assume readers have background in ordinary differential equations, multivariable calculus, and basic stochastic calculus, as well as some exposure to numerical methods for SDEs~\cite{higham2021, kloeden1992, oksendal2003}. Familiarity with differentiation of ODE-based objectives (the discretize-optimize and optimize-discretize approaches) is helpful but not required. Depending on the course and students' backgrounds, instructors may consider teaching the differentiation through ODEs before presenting the stochastic case; see~\cite{kirk2004} for background on these techniques in the deterministic setting. 

To keep the material accessible, we have opted for a mostly informal presentation style and small-dimensional numerical examples, referring to rigorous results and more complex examples where suitable. Our presentation draws substantially from the comprehensive treatment in Kidger~\cite{kidger2022neural}, which provides an elegant and abstract framework using rough path theory and controlled differential equations to unify the treatment of ODEs and SDEs. We have adapted this material into a more simplified approach that prioritizes accessibility and intuitive understanding. 

The paper is organized into two progressive tales, corresponding to the two classic differentiation approaches for SDEs. Tale 1 covers the discretize-optimize approach, showing how to differentiate through numerical schemes for both It\^{o} and Stratonovich SDEs --- the two standard frameworks for stochastic calculus. Tale 2 examines the optimize-discretize approach, investigating whether the equations from Tale 1 converge to well-defined continuous processes. Throughout the paper, we validate our theoretical results with numerical experiments using the Black--Scholes model. The code for reproducing the figures in this paper and validating the theoretical results is available at \url{https://github.com/EmoryMLIP/2026-DifferentiatingSDEPrimer}.

The material naturally fits into courses on stochastic processes, quantitative finance, scientific computing, and machine learning. While standard courses teach SDE simulation, they rarely cover gradient computation despite its importance in applications. This tutorial fills that gap by connecting simulation methods students already know to the techniques needed for optimization and sensitivity analysis. Instructors can use the individual tales as standalone modules or present both tales together for comprehensive coverage.

By the end of this primer, readers will understand how to differentiate through SDEs correctly in both It\^{o} and Stratonovich settings using either discretize-optimize or optimize-discretize. We seek to teach readers to recognize why endpoint evaluation matters and understand the fundamental differences between differentiating deterministic and stochastic systems. This foundation enables reasoning about stochastic differentiation problems beyond the cases we explicitly cover.

\section{Tale 1. Differentiating through SDEs. Discretize-Optimize} We demonstrate the reliability of the discretize-optimize approach for computing gradients through SDEs. We show how to discretize It\^{o} and Stratonovich SDEs using corresponding standard numerical schemes, then differentiate through their chosen respective schemes via automatic differentiation to obtain the desired sensitivities.

We adopt a notational convention where boldface is used for matrices to distinguish them from random variables, both of which appear as uppercase letters, and vectors are represented by lowercase letters. This convention naturally helps maintain visual distinction between matrices, vectors, and stochastic processes.

\subsection{Setup and Objective}

Suppose we have an objective function that depends on $X(T)$, the state of a stochastic system at time $T$ starting from the initial state $x_0 \in \mathbb{R}^d$,
\begin{equation}\label{eq:objective_SDE_Ito}
    J(x_0) = \mathbb{E}\left[ \Phi(X(T)) \ | \ X(0) = x_0 \right],
\end{equation}
where $\Phi : \mathbb{R}^d \to \mathbb{R}$ is a smooth function. Our task is to compute the gradient of this objective function with respect to the initial state $x_0$. As most stochastic systems do not admit explicit solutions, one has to rely virtually always on numerical approximations of the gradients. 

There are two natural ways to proceed. In the first approach, one discretizes the state dynamics and directly differentiates the corresponding discrete objective function. In the second approach, one first derives a continuous adjoint differential equation whose initial state provides the gradient of the objective function, then discretizes that equation and numerically integrates. The first approach is called \textit{discretize-optimize}, whereas the second approach is called \textit{optimize-discretize} (known in the PDE-constrained optimization literature as \textit{discretize-then-optimize} and \textit{optimize-then-discretize}, respectively). In this tale, we tell the story of discretize-optimize and show that it is reliable for the two standard frameworks for stochastic systems: It\^{o} and Stratonovich. 

\subsection{It\^{o} SDEs. Forward and Reverse-Mode Differentiation}

Assume that the state dynamics evolves according to the It\^{o} SDE
\begin{equation}\label{eq:SDE_Ito}
    dX(s) = f(s,X(s)) ds + g(s,X(s))dW(s), \quad 0 \leq s \leq T, \quad X(0) = x_0,
\end{equation}
where $f : [0,T] \times \mathbb{R}^d \to \mathbb{R}^d$ is the drift, $g : [0,T] \times \mathbb{R}^d \to \mathbb{R}^{d \times m}$ is the diffusion, and $W(s)$ is an $m-$dimensional Brownian motion. Together, these terms describe a system that evolves under both deterministic and stochastic effects. 

En route to differentiating the objective~\eqref{eq:objective_SDE_Ito} via discretize-optimize, we discretize the forward dynamics first, then differentiate the resulting discrete objective function. Using the Euler-Maruyama scheme~\cite[Chap.~8]{higham2021}, \cite[Sec.~10.2]{kloeden1992} with timestep $\Delta t = T/N$ and increments $\Delta W_n \sim \mathcal{N}(0, \Delta t \mathbf{I}_m)$ we discretize the It\^{o} SDE via
\begin{equation}\label{eq:SDE_discrete_Ito}
    X_{n+1} = X_n + \Delta t f(t_n,X_n) +  g(t_n,X_n) \Delta W_n, \quad 0\leq n \leq N-1,\quad X_0=x_0.
\end{equation}
Thus, the discrete objective function is given by
\begin{equation*}\label{eq:objective_SDE_Ito_discrete}
    J_{\Delta t}(x_0) = \mathbb{E}[\Phi(X_N) \ | \ X_0 = x_0].
\end{equation*}

Under mild regularity conditions the interchange of differentiation and expectation holds, so the correct gradient of our continuous objective is an average of pathwise gradients, where $X_0 = x_0$ is the deterministic initial state, as follows
\begin{equation*}
    \nabla_{x_0} J_{\Delta t}(x_0) = \mathbb{E}[\nabla_{x_0} \Phi(X_N) \ | \ X_0 = x_0].
\end{equation*}

Hence, our task reduces to computing the pathwise gradient $\nabla_{x_0} \Phi(X_N)$ over many simulated paths, and since $\Phi(X_N)$ is a finite composition, the chain rule applies. Using \textit{forward-mode differentiation}, we could apply the chain rule directly by first computing the Jacobian
\begin{equation}\label{eq:X_N_X_0_jac_discrete}
    \frac{\partial X_N}{\partial X_0}=\frac{\partial X_N}{\partial X_{N-1}}~\frac{\partial X_{N-1}}{\partial X_{N-2}} ~\cdots ~\frac{\partial X_1}{\partial x_{0}},
\end{equation}
where differentiating the Euler-Maruyama scheme~\eqref{eq:SDE_discrete_Ito} with respect to $X_n$ yields the one-step Jacobian
\begin{equation}\label{eq:one_step_jacobian_euler}
    \frac{\partial X_{n+1}}{\partial X_{n}} = \mathbf{I} + \Delta t \frac{\partial f(t_n,X_n)}{\partial X_n} +  \Delta W_n \frac{\partial g(t_n,X_n)}{\partial X_n},
\end{equation}
which we can plug into
\begin{equation}\label{eq:J_grad_discrete_direct}
    \nabla_{x_0} J_{\Delta t}(x_0) = \mathbb{E}[\nabla_{x_0} \Phi(X_N) \ | \ X_0 = x_0] = \mathbb{E} \left[ \frac{\partial \Phi(X_N)}{\partial X_N} \frac{\partial X_N}{\partial x_0} \  \bigg| \ X_0 = x_0 \right].
\end{equation}
Alternatively, after simulating trajectories, one could use \textit{reverse-mode automatic differentiation (backpropagation)}~\cite[Chap.~3]{griewank2008} which reorganizes the chain rule into a backward recursion over the time steps. More precisely, instead of computing full Jacobian matrices forward in time, we compute a matrix-vector product with its transpose backward in time. To this end, we define the \textit{discrete adjoint} $p_n \in \mathbb{R}^d$ as the gradient of the discrete pathwise objective with respect to the state at step $n$; that is,
\begin{equation*}
p_n=\frac{\partial \Phi(X_N)}{\partial X_{n}}^\top =\frac{\partial X_{n+1}}{\partial X_{n}}^\top \frac{\partial \Phi(X_N)}{\partial X_{n+1}}^\top=\frac{\partial X_{n+1}}{\partial X_{n}}^\top p_{n+1}.
\end{equation*}
Hence, using~\eqref{eq:one_step_jacobian_euler} we obtain the discrete adjoint recursion
\begin{equation} \label{eq:ito_discrete_adj}
    \begin{cases}
        p_n = \left( \mathbf{I} + \Delta t \frac{\partial f(t_n,X_n)}{\partial X_n}^\top + \Delta W_n^\top \frac{\partial g(t_n,X_n)}{\partial X_n}^\top \right)  p_{n+1},\quad 0\leq n \leq N-1 \\ 
        p_N = \frac{\partial \Phi(X_N)}{\partial X_{N}}^\top,
    \end{cases}
\end{equation}
where the gradient is thus given by
\begin{equation} \label{eq:grad_j_discrete_reverse}
    \nabla_{x_0} J_{\Delta t}(x_0)= \mathbb{E}[p_0 \ | \ X_0 = x_0].
\end{equation}

Mathematically, both formulations~\eqref{eq:X_N_X_0_jac_discrete},~\eqref{eq:J_grad_discrete_direct} and~\eqref{eq:ito_discrete_adj},~\eqref{eq:grad_j_discrete_reverse} yield the same gradient; once the discrete objective function is defined, either method computes it.

The key advantage of~\eqref{eq:ito_discrete_adj},~\eqref{eq:grad_j_discrete_reverse} is computational. The direct chain-rule expansion~\eqref{eq:X_N_X_0_jac_discrete} propagates the full sensitivity matrices $\frac{\partial X_{n}}{\partial x_{0}} \in \mathbb{R}^{d \times d}$, requiring repeated Jacobian multiplications at each time step. The backward sweep, by contrast, updates a single vector $p_n \in \mathbb{R}^d$ using matrix–vector products $\Bigl(\tfrac{\partial X_{n+1}}{\partial X_n}\Bigr)^{\top} p_{n+1}$, reducing the per-step computational cost considerably. Hence, reverse-mode automatic differentiation  organizes the chain rule in the most computationally efficient way.

\begin{summarybox}[Summary: Discrete Adjoint for It\^{o} SDEs (Euler-Maruyama)] Assume that the It\^{o} state dynamics is discretized via Euler-Maruyama scheme
    \[
    X_{n+1} = X_n + \Delta t f(t_n,X_n) +  g(t_n,X_n) \Delta W_n, \quad 0\leq n \leq N-1,\quad X_0=x_0.
    \]
    Then the gradient of the discrete objective function
    \[
    J_{\Delta t}(x_0) = \mathbb{E}[\Phi(X_N) \ | \ X_0 = x_0]
    \]
    can be computed via reverse-mode automatic differentiation as
    \[
    \nabla_{x_0} J_{\Delta t}(x_0)= \mathbb{E}[p_0 \ | \ X_0 = x_0],
    \]
    where the pathwise discrete adjoint state, $p_n$, satisfies the backward recursion
    \[
    \begin{cases}
        p_n = \left( \mathbf{I} + \Delta t \frac{\partial f(t_n,X_n)}{\partial X_n}^\top + \Delta W_n^\top \frac{\partial g(t_n,X_n)}{\partial X_n}^\top \right)  p_{n+1},\quad 0\leq n \leq N-1 \\ 
        p_N = \frac{\partial \Phi(X_N)}{\partial X_{N}}^\top.
    \end{cases}
    \]
\end{summarybox}

\subsection{Parameters and Running Costs}\label{subsec:params_running} So far, we have focused on computing gradients with respect to the initial state $x_0$. However, in many applications, we are also interested in sensitivities with respect to parameters that appear in the SDE itself. One example would be a neural SDE where the goal is to optimize the weights of the neural network. Another example is in finance, where one might want to compute sensitivities such as Vega $(\partial C/\partial \sigma)$ or Rho $(\partial C/\partial r)$; we will see this in our Black--Scholes example in \hyperref[tab:greeks]{Table~\ref*{tab:greeks}}.

A simple and effective way to derive sensitivities for parameters is to add them as states.  Consider an It\^{o} SDE with parameters $\theta \in \mathbb{R}^p$:
\begin{equation*}
    dX(s) = f(s,X(s),\theta) ds + g(s,X(s),\theta)dW(s), \quad 0 \leq s \leq T, \ X(0) =x_0
\end{equation*}
with objective
\begin{equation} \label{eq:param_objective}
    J(x_0,\theta) = \mathbb{E}\left[ \Phi(X(T; x_0, \theta)) \right].
\end{equation}

To compute gradients of~\eqref{eq:param_objective} with respect to both $x_0$ and $\theta$, we augment the state; that is, we treat $\theta$ as an additional state variable that remains constant in time. Using Euler-Maruyama as before, we obtain
\begin{equation*}
    \begin{cases}
        X_{n+1} = X_n + \Delta t f(t_n,X_n,\theta_n) +  g(t_n,X_n,\theta_n) \Delta W_n \\
        \theta_{n+1} = \theta_n, \quad 0\leq n \leq N-1,\quad \theta_0=\theta.
    \end{cases}
\end{equation*}
To derive the discrete adjoint for this augmented system, we define the augmented state $\widehat{X}_n=(X_n, \theta_n)^\top$ and apply the chain rule as before,
\begin{equation} \label{eq:chain_rule_parameters}
    \frac{\partial \Phi(X_N)}{\partial \widehat{X}_n}^\top = \left(\frac{\partial \widehat{X}_{n+1}}{\partial \widehat{X}_n}\right)^\top \frac{\partial \Phi(X_N)}{\partial \widehat{X}_{n+1}}^\top.
\end{equation}
The Jacobian of the augmented update is a block matrix:
\begin{equation*}
    \frac{\partial \widehat{X}_{n+1}}{\partial \widehat{X}_n} = \begin{pmatrix}
        \frac{\partial X_{n+1}}{\partial X_n} & \frac{\partial X_{n+1}}{\partial \theta_n} \\
        \mathbf{0} & \mathbf{I}
    \end{pmatrix},
\end{equation*}
where the top-left block is the one-step Jacobian from~\eqref{eq:one_step_jacobian_euler}, the top-right block captures how $X_{n+1}$ depends on $\theta$, the bottom-left block is zero because $\theta_{n+1}$ does not depend on $X_n$, and the bottom-right block is the identity because $\theta_{n+1} = \theta_n$. Now,  defining $p_n, q_n$ as the discrete adjoints for state and parameters respectively
\begin{equation*}
    p_n = \frac{\partial \Phi(X_N)}{\partial X_n}^\top, \qquad q_n = \frac{\partial \Phi(X_N)}{\partial \theta_n}^\top,
\end{equation*}
substituting into~\eqref{eq:chain_rule_parameters} yields:
\begin{equation*}
    \begin{pmatrix} p_n \\ q_n \end{pmatrix} = \begin{pmatrix}
        \frac{\partial X_{n+1}}{\partial X_n}^\top & 0 \\
        \frac{\partial X_{n+1}}{\partial \theta_n}^\top & \mathbf{I}
    \end{pmatrix} \begin{pmatrix} p_{n+1} \\ q_{n+1} \end{pmatrix}.
\end{equation*}
Multiplying this out gives a discrete adjoint with two components, one for the state and the other for the parameters:
\begin{equation} \label{eq:discrete_adjoint_param}
    \begin{cases}
        p_n = \left( \mathbf{I} + \Delta t \frac{\partial f(t_n,X_n,\theta)}{\partial X_n}^\top + \Delta W_n^\top \frac{\partial g(t_n,X_n,\theta)}{\partial X_n}^\top \right)  p_{n+1} \\
        q_n = q_{n+1} + \Delta t \frac{\partial f(t_n,X_n,\theta)}{\partial \theta}^\top p_{n+1} + \Delta W_n^\top \frac{\partial g(t_n,X_n,\theta)}{\partial \theta}^\top p_{n+1} \\
        p_N = \frac{\partial \Phi(X_N)}{\partial X_{N}}^\top, \quad q_N = 0,
    \end{cases}
\end{equation}
where $q_N = 0$ since the terminal cost $\Phi$ does not directly depend on the parameters $\theta$. More generally, if $\Phi$ depends explicitly on $\theta$, then $q_N = \frac{\partial \Phi(X_N, \theta)}{\partial \theta}^\top$. The gradients are then given by
\begin{equation} \label{eq:param_gradients}
\nabla_{x_0} J_{\Delta t}(x_0,\theta)= \mathbb{E}[p_0 \ | \ X_0 = x_0], \quad \nabla_{\theta} J_{\Delta t}(x_0,\theta)= \mathbb{E}[q_0 \ | \ X_0 = x_0].
\end{equation}

The key insight is that in~\eqref{eq:discrete_adjoint_param},~\eqref{eq:param_gradients} parameter gradients accumulate information from all time steps through the $q_n$ recursion, while state gradients $p_n$ simultaneously propagate backward through the chain rule as usual. Automatic differentiation frameworks can handle this bookkeeping automatically, computing both $p_0$ and $q_0$ in a single backward pass.

The same strategy can also be used when the objective function includes running costs, which depend on the entire trajectory, not just the terminal state. Such objectives take the form
$$
J(x_0) = \mathbb{E}\left[\Phi(X(T; x_0)) + \int_0^T L(s, X(s)) ds\right],
$$
where $L(s, X(s))$ represents an instantaneous contribution and is commonly referred to as the \textit{running cost}.
To handle this, we can  introduce an auxiliary state variable $Y(s)$ that accumulates the running cost along trajectories:
$$
dY(s) = L(s, X(s)) ds, \quad  \quad 0 \leq s \leq T, \ Y(0) = 0,
$$
with terminal value $Y(T; x_0) = \int_0^T L(s, X(s)) ds$. Augmenting the state to $\widehat{X}(s) = (X(s), Y(s))^\top$ reduces the problem to the terminal cost case with objective
$$
J(x_0) =\mathbb{E}[\Phi(X(T; x_0)) + Y(T; x_0)] = \mathbb{E}[\widehat{\Phi}(\widehat{X}(T; x_0))],
$$
where $\widehat{\Phi}(\widehat{X}(T; x_0)) = \Phi(X(T; x_0)) + Y(T; x_0)$. This is now a terminal cost problem for the augmented system, and the adjoint techniques from this tale apply directly.

\subsection{Example. Black--Scholes SDE}

To illustrate the above derivative computation, we use the Black--Scholes model~\cite{black1973,merton1973} as our primary example because it is often covered in the curriculum and derivatives are important for option pricing. This model describes the evolution of a stock price $S(t)$ under geometric Brownian motion as an It\^{o} SDE
\begin{equation}\label{eq:BS_Ito}
    dS(t) = \mu S(t) dt + \sigma S(t) dW(t), \quad S(0) = S_0,
\end{equation}
where $\mu$ is the drift (expected return), $\sigma > 0$ is the volatility, and $S_0 > 0$ is the initial stock price. This SDE admits an analytical solution~\cite{black1973}
\begin{equation*}
S(t) = S_0 \exp\left(\left(\mu - \frac{\sigma^2}{2}\right)t + \sigma W(t)\right),
\end{equation*}
which will prove useful for validating our numerical computations.
In practice, the expected returns and volatility are usually time-dependent, which means there is no analytical solution and numerical simulation and differentiation are required.

In options pricing, one seeks not only the fair value of an option but also its sensitivities to market parameters, quantities known collectively as the Greeks~\cite{hull2018}. For instance, consider a European call option with strike price $K$ and maturity $T$. Under the risk-neutral measure (where $\mu = r$, the risk-free rate), the option value at time $t=T$ is given by
\begin{equation} \label{eq:option_pricing}
    C(S_0, K, T) = \mathbb{E}\left[\max(S(T) - K, 0)\right].
\end{equation}

The Greeks quantify how this option value responds to changes in the underlying parameters. The most frequently computed Greeks are listed in \hyperref[tab:greeks]{Table~\ref*{tab:greeks}}.

\begin{table}[ht]
    \centering
    \begin{tabular}{lll}
    \toprule
    Greek & Definition & Interpretation \\
    \midrule
    Delta & $\Delta = \partial C/\partial S_0$ & Sensitivity to stock price \\
    Gamma & $\Gamma = \partial^2 C/\partial S_0^2$ & Curvature with respect to stock price \\
    Vega & $\mathcal{V} = \partial C/\partial \sigma$ & Sensitivity to volatility \\
    Theta & $\Theta = -\partial C/\partial T$ & Time decay \\
    Rho & $\rho = \partial C/\partial r$ & Sensitivity to interest rate \\
    \bottomrule
    \end{tabular}
    \caption{Common Greeks in option pricing and their financial interpretations}
    \label{tab:greeks}
\end{table}

\FloatBarrier

As computing the Greeks is precisely a problem of differentiating an objective with respect to parameters that appear in the underlying SDE~\eqref{eq:BS_Ito}, this makes it the appropriate vehicle for illustrating the discretize-optimize and optimize-discretize approaches. We will return to this example in the sections that follow.

To validate our derivative computations thus far, we conduct a numerical experiment using option pricing as our test case. More specifically, we compute the Greek $\Delta = \partial C / \partial S_0 $, where the option price is given by \eqref{eq:option_pricing}. Using the parameters $S_0 = 100$, $K = 110$, $r = 0.05$, $\sigma = 0.2$, and $T = 1$, we simulate $10^5$ Monte Carlo paths with varying time step sizes $\Delta t \in [10^{-3}, 1]$. For each $\Delta t$, we compute the discrete adjoint gradient, average over all paths, and compare against the known analytical Black--Scholes $\Delta$ obtained from the closed-form formula~\cite{black1973}.

\begin{figure}[htp] 
  \centering
  \includegraphics[width=.78\linewidth]{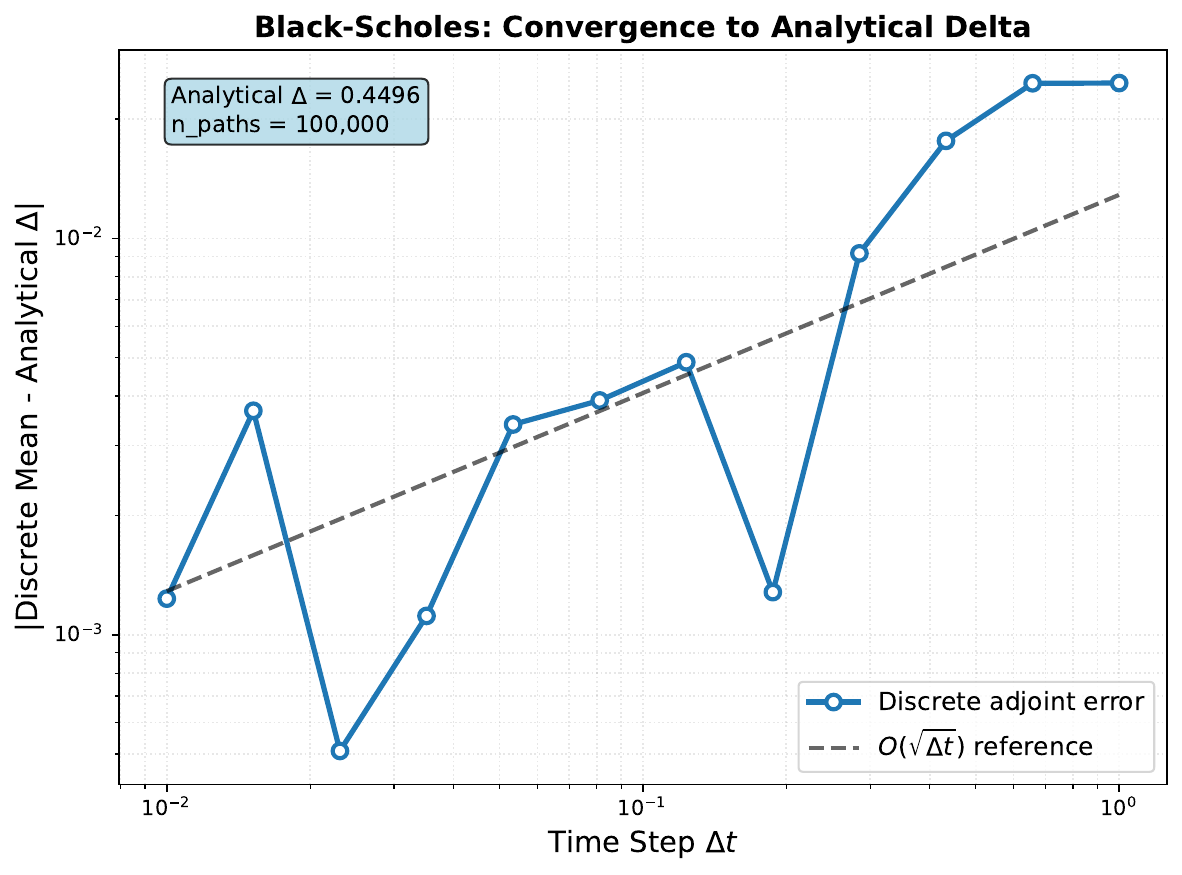}
  \caption{Error in discrete adjoint gradient compared to the analytical $\Delta$ versus time step size $\Delta t$ for Black--Scholes model on log-log scale. Computations using $n = 10^5$ Monte Carlo paths with parameters $S_0 = 100$, $K = 110$, $r = 0.05$, $\sigma = 0.2$, $T = 1$.}
  \label{fig:ito_convergence}
\end{figure}

\hyperref[fig:ito_convergence]{Figure~\ref*{fig:ito_convergence}} shows the convergence behavior. As $\Delta t$ decreases, the error between the approximation and the true Greek $\Delta$, $|\overline{\Delta}_{\text{discrete}} - \Delta_{\text{analytical}}|$, diminishes following the expected $O(\sqrt{\Delta t})$ rate characteristic of the Euler-Maruyama scheme~\cite[Chap.~8.4]{higham2021}, \cite[Thm.~10.2.2]{kloeden1992}. For small $\Delta t < 0.01$, the errors plateau around $10^{-3}$ due to Monte Carlo sampling noise rather than discretization error.

\subsection{Stratonovich SDEs}

Having illustrated that discretize-optimize produces reliable gradients for It\^{o} SDEs (Euler-Maruyama), we now turn to the Stratonovich framework. Recall that the key distinction between It\^{o} and Stratonovich calculus stems from the definition of the corresponding stochastic integrals. More specifically, stochastic Riemann sums use different evaluation points of the integrands: It\^{o} uses a left-endpoint evaluation (Euler-Maruyama), while Stratonovich uses a midpoint evaluation. A consequence is that Stratonovich SDEs satisfy the classical chain rule; see~\cite[Chap.~7]{higham2021} for a detailed comparison. Assume that we have the same objective function~\eqref{eq:objective_SDE_Ito}
but now the state $X(s)$ evolves according to the Stratonovich SDE
\begin{equation}\label{eq:SDE_Stratonovich}
dX(s) = f(s,X(s)) ds + g(s,X(s)) \circ dW(s), \quad 0 \leq s \leq T, \ X(0) = x_0.
\end{equation}

To enact the discretize-optimize approach we need to discretize this SDE, but the Euler-Maruyama scheme employed before is no longer suitable as it converges to an It\^{o} SDE~\cite[Chap.~8]{higham2021}, \cite[Sec.~10.2]{kloeden1992}. To recover the Stratonovich SDE in the limit, we use the Heun scheme~\cite[Chap.~17.3]{higham2021} with timestep $\Delta t = T/N$ and Brownian increments $\Delta W_n \sim \mathcal{N}(0, \Delta t \mathbf{I}_m)$ as follows
\begin{equation*}
    \begin{split}
        X_{n+1} &= X_n + \frac{\Delta t}{2} \left[ f(t_n,X_n) + f(t_{n+1},\widetilde{X}_{n+1}) \right] +  \left[ g(t_n,X_n) + g(t_{n+1},\widetilde{X}_{n+1}) \right] \frac{\Delta W_n}{2} \\
        \widetilde{X}_{n+1} &= X_n + \Delta t  f(t_n,X_n) + g(t_n,X_n) \Delta W_n. 
    \end{split}
\end{equation*}

As before, the discrete objective function is
\begin{equation*}
    J_{\Delta t}(x_0)=\mathbb{E}\left[\Phi(X_N) \mid X_0=x_0\right],
\end{equation*}
which corresponds to pathwise objective function $\Phi(X_N)$. To compute the pathwise gradients, we introduce the discrete adjoint $p_n \in \mathbb{R}^d$ and derive the corresponding reverse-mode automatic differentiation scheme.

Note that the Heun scheme involves an intermediate predictor step $\widetilde{X}_{n+1}$, so we must carefully apply the chain rule through this two-stage update. Differentiating the predictor step with respect to $X_n$ gives
\begin{equation*}
    \frac{\partial \widetilde{X}_{n+1}}{\partial X_n} = \mathbf{I} + \Delta t\frac{\partial f(t_n,X_n)}{\partial X_n} + \frac{\partial g(t_n,X_n)}{\partial X_n} \Delta W_n.
\end{equation*}
For the corrector step, differentiating $X_{n+1}$ with respect to $X_n$ requires the chain rule through $\widetilde{X}_{n+1}$; that is,
\begin{align*}
    \frac{\partial X_{n+1}}{\partial X_n} &= \mathbf{I} + \frac{\Delta t}{2} \frac{\partial f(t_n,X_n)}{\partial X_n} + \frac{\partial g(t_n,X_n)}{\partial X_n} \frac{\Delta W_n}{2} \\
    & \ +  \frac{\Delta t}{2} \frac{\partial f(t_{n+1},\widetilde{X}_{n+1})}{\partial \widetilde{X}_{n+1}} \frac{\partial \widetilde{X}_{n+1}}{\partial X_n} + \frac{\partial g(t_{n+1},\widetilde{X}_{n+1})}{\partial \widetilde{X}_{n+1}} \frac{\partial \widetilde{X}_{n+1}}{\partial X_n} \frac{\Delta W_n}{2}.
\end{align*}

Substituting the predictor Jacobian and collecting terms, we obtain
\begin{equation*}
\begin{split}
    \frac{\partial X_{n+1}}{\partial X_n} = & \ \mathbf{I} + \frac{\Delta t}{2} \left[\frac{\partial f(t_n,X_n)}{\partial X_n} + \frac{\partial f(t_{n+1},\widetilde{X}_{n+1})}{\partial \widetilde{X}_{n+1}}\right]\\
    &+  \left[\frac{\partial g(t_n,X_n)}{\partial X_n} + \frac{\partial g(t_{n+1},\widetilde{X}_{n+1})}{\partial \widetilde{X}_{n+1}}\right] \frac{\Delta W_n}{2}+H_n,
\end{split}
\end{equation*}
where $H_n$ collects the additional terms from the predictor-corrector structure
\begin{equation} \label{eq:higher_order_terms_heun}
\begin{split}
    H_n = {}&\frac{\Delta t^2}{2} \frac{\partial f(t_{n+1},\widetilde{X}_{n+1})}{\partial \widetilde{X}_{n+1}} \frac{\partial f(t_n,X_n)}{\partial X_n}\\
    &+ \frac{(\Delta W_n)^2}{2} \frac{\partial g(t_{n+1},\widetilde{X}_{n+1})}{\partial \widetilde{X}_{n+1}} \frac{\partial g(t_n,X_n)}{\partial X_n}\\
    &+ \frac{\Delta t \Delta W_n}{2} \frac{\partial f(t_{n+1},\widetilde{X}_{n+1})}{\partial \widetilde{X}_{n+1}} \frac{\partial g(t_n,X_n)}{\partial X_n}\\
    &+ \frac{\Delta t \Delta W_n}{2} \frac{\partial g(t_{n+1},\widetilde{X}_{n+1})}{\partial \widetilde{X}_{n+1}} \frac{\partial f(t_n,X_n)}{\partial X_n}.
\end{split}
\end{equation}
Thus, the discrete adjoint recursion is
\begin{equation*}
    p_n = \left(\frac{\partial X_{n+1}}{\partial X_n}\right)^\top p_{n+1}, 
\end{equation*}
which, using the Jacobian above, becomes
\begin{equation}\label{eq:heun_adjoint}
    \begin{cases}
        p_n = \bigg( \mathbf{I} + \frac{\Delta t}{2} \left[\frac{\partial f(t_n,X_n)}{\partial X_n}^\top + \frac{\partial f(t_{n+1},\widetilde{X}_{n+1})}{\partial \widetilde{X}_{n+1}}^\top\right]\\ \qquad
        +  \frac{\Delta W_n^\top}{2} \left[\frac{\partial g(t_n,X_n)}{\partial X_n}^\top + \frac{\partial g(t_{n+1},\widetilde{X}_{n+1})}{\partial \widetilde{X}_{n+1}}^\top\right] + H_n^\top \bigg) p_{n+1},\  0\leq n \leq N-1 \\
        p_N = \frac{\partial \Phi(X_N)}{\partial X_N}^\top
    \end{cases}
\end{equation}
The gradient is hence given by
\begin{equation*}
    \nabla_{x_0} J_{\Delta t}(x_0)=\mathbb{E}\left[p_0 \mid X_0=x_0\right].
\end{equation*}
Note that while the discrete adjoint for Stratonovich SDEs involves more terms than its It\^{o} counterpart due to the predictor-corrector structure, modern automatic differentiation frameworks handle this complexity automatically, computing the full expression~\eqref{eq:heun_adjoint} including the higher order $\Delta t^2$ and $\Delta t \Delta W_n$ terms in $H_n$. These higher-order terms cannot be dropped without changing what is computed: they are part of the chain rule applied to the Heun update, so keeping them is what makes $p_0$ the exact gradient of the discrete objective rather than the gradient of a simplified scheme.

The same augmentation strategy from~\Cref{subsec:params_running} applies to Stratonovich SDEs discretized with the Heun scheme since the augmented process $(X,\theta)$ remains a  Stratonovich SDE. The parameter adjoint $q_n$ accumulates gradient information from all time steps in exactly the same way, with the only difference being the one-step Jacobian used in the $p_n$ recursion. 

\begin{summarybox}[Summary: Discrete Adjoint for Stratonovich SDEs (Heun)]
Assume that the Stratonovich state dynamics is discretized via the Heun scheme
    \begin{equation*}
    \begin{split}
        \widetilde{X}_{n+1} &= X_n + \Delta t  f(t_n,X_n) + g(t_n,X_n) \Delta W_n \\
        X_{n+1} &= X_n + \frac{\Delta t}{2} \left[ f(t_n,X_n) + f(t_{n+1},\widetilde{X}_{n+1}) \right] \\
        &+  \left[ g(t_n,X_n) + g(t_{n+1},\widetilde{X}_{n+1}) \right] \frac{\Delta W_n}{2}.
    \end{split}
    \end{equation*}
    Then the gradient of the discrete objective function
    \[
    J_{\Delta t}(x_0) = \mathbb{E}[\Phi(X_N) \ | \ X_0 = x_0]
    \]
    can be computed via reverse-mode automatic differentiation as
    \[
    \nabla_{x_0} J_{\Delta t}(x_0)= \mathbb{E}[p_0 \ | \ X_0 = x_0],
    \]
    where the pathwise discrete adjoint state, $p_n$, satisfies the backward recursion
    \[
    \begin{cases}
        p_n =\bigg( \mathbf{I} + \frac{\Delta t}{2} \left[\frac{\partial f(t_n,X_n)}{\partial X_n}^\top + \frac{\partial f(t_{n+1},\widetilde{X}_{n+1})}{\partial \widetilde{X}_{n+1}}^\top\right]\\ \qquad
    +  \frac{\Delta W_n^\top}{2} \left[\frac{\partial g(t_n,X_n)}{\partial X_n}^\top + \frac{\partial g(t_{n+1},\widetilde{X}_{n+1})}{\partial \widetilde{X}_{n+1}}^\top\right] + H_n^\top \bigg) p_{n+1},\quad 0\leq n \leq N-1 \\
        p_N = \frac{\partial \Phi(X_N)}{\partial X_N}^\top.
    \end{cases}
    \]
    with $H_n$ defined in~\eqref{eq:higher_order_terms_heun}.
\end{summarybox} 

\subsection{The Lesson. Discretize-Optimize is Safe for SDEs} 

While stochastic dynamical systems require careful selection of the discretization scheme for the forward dynamics, the practical takeaway is reassuring. For SDEs, differentiating directly through the numerical scheme via backpropagation is not an approximation or shortcut; it is the most direct way to compute sensitivities of an objective. The standard ODE workflow --- simulate forward, then backpropagate --- applies seamlessly to the stochastic case, except that gradients must be averaged over random trajectories.

\section{Tale 2. Use Stratonovich for Optimize-Discretize}\label{sec:tale_2}

So far, we have derived discrete adjoints for both It\^{o} and Stratonovich SDEs by differentiating through their respective numerical schemes and learned that discretize-optimize is a safe choice for differentiating through SDEs. But a natural mathematical question arises: as the timestep $\Delta t \to 0$, do these discrete adjoints converge to well-defined continuous adjoint processes? And if so, what are those continuous processes?

The answers to these questions are important for at least two reasons. First, in essence, the existence of a well-defined continuous process tells us that our computational framework is stable under mesh refinement. This means that things \enquote{will not break} and have a well-defined limit if we consider a more accurate discretization of the state dynamics. Second, the existence of a continuous process presents alternative approximation opportunities for the discrete adjoint. The latter is the core of the optimize-discretize approach: one first derives a continuous adjoint and then discretizes it.

\subsection{From Discrete to Continuous Adjoints} To build intuition, we briefly review the deterministic case. Taking $g \equiv 0$ in our It\^{o}~\eqref{eq:SDE_Ito} and Stratonovich~\eqref{eq:SDE_Stratonovich} SDEs, they both reduce to the ODE
\begin{equation*}
    \frac{dx(s)}{ds} = f(s,x(s)), \quad s\in (0,T),\quad  x(0) = x_0.
\end{equation*}
The expectation in the objective~\eqref{eq:objective_SDE_Ito} becomes trivial since the system is deterministic, reducing to 
\begin{equation}\label{eq:objective_ODE}
    J(x_0) = \Phi(x(T; x_0)).
\end{equation}
Thus, if we want to compute the gradient of this objective function with respect to the initial state in this setting, the discrete adjoint~\eqref{eq:ito_discrete_adj} from It\^{o} (Euler-Maruyama) simplifies to the backward discrete adjoint recursion corresponding to the objective function~\eqref{eq:objective_ODE} discretized with forward Euler
\begin{equation} \label{eq:ODEadj_rec_d-o_euler}
    \begin{cases}
        p_n = \left(\mathbf{I}+\Delta t\frac{\partial f(t_{n},x_{n})}{\partial x_{n}}^\top\right) p_{n+1},\quad 0\leq n \leq N-1\\ 
        p_N = \frac{\partial \Phi(x_N)}{\partial x_N}^\top,
    \end{cases}
\end{equation}
while the discrete adjoint~\eqref{eq:heun_adjoint} from Stratonovich (Heun) to the backward recursion corresponding to the discrete objective function discretized with modified Euler
\begin{equation}\label{eq:ODEadj_rec_d-o_mod_euler}
    \begin{cases}
        p_n =\bigg( \mathbf{I} + \frac{\Delta t}{2} \left[\frac{\partial f(t_n,x_n)}{\partial x_n}^\top + \frac{\partial f(t_{n+1},\widetilde{x}_{n+1})}{\partial \widetilde{x}_{n+1}}^\top\right] + D_n^\top \bigg) p_{n+1}\\

        p_N = \frac{\partial \Phi(x_N)}{\partial x_N}^\top,
    \end{cases}
\end{equation}
where $D_n$ contains the deterministic version of the predictor-corrector correction term
\begin{equation*}
    D_n = \frac{\Delta t^2}{2}\frac{\partial f(t_{n+1},\widetilde{x}_{n+1})}{\partial \widetilde{x}_{n+1}} \frac{\partial f(t_n,x_n)}{\partial x_n}.
\end{equation*}

Both recursions can be interpreted as Riemann-sum-like approximations (quadrature rules) when unrolled. The two schemes above, however, are not standard Riemann sums (left, right, or midpoint):~\eqref{eq:ODEadj_rec_d-o_euler} pairs the Jacobian at $t_n$ with the adjoint at $t_{n+1}$, while~\eqref{eq:ODEadj_rec_d-o_mod_euler} involves a trapezoidal average of Jacobians at $t_n$ and $t_{n+1}$, together with higher-order correction terms, all multiplied by the adjoint at $t_{n+1}$. For sufficiently smooth ODE trajectories, both are nevertheless consistent discretizations and converge to the same continuous adjoint ODE~\cite[Thm.~2.1]{hager2000}, via the transformed adjoint system, and~\cite[Thm.~5.1]{walther2007}, via reverse-mode differentiation of the Runge-Kutta scheme
\begin{equation}\label{eq:p_ODE_limit}
    \begin{cases}
        \dot{p}(s)+ \frac{\partial f(s,x(s))}{\partial x}^\top p(s)=0,\quad s\in (0,T)\\
        p(T)= \frac{\partial \Phi (x(T))}{\partial x}^\top,
    \end{cases}
\end{equation}
where the variable $p$ in~\eqref{eq:p_ODE_limit} is called the \textit{continuous adjoint state}~\cite[Chap.~5.1]{kirk2004} and our desired gradient is given by
\begin{equation*}
    \nabla_{x_0} J(x_0) = p(0).
\end{equation*}
The convergence rate inherits the order of the forward Runge-Kutta scheme~\cite[Thm.~5.1]{walther2007}; in our case, the positive Runge-Kutta weights of the forward Euler and Heun schemes preserve this order.

This convergence is precisely what makes optimize-discretize work for ODEs: one can derive the continuous adjoint~\eqref{eq:p_ODE_limit} directly, discretize it in any consistent way, and the resulting numerical approximation will converge to $p(0)$ as $\Delta t \to 0$. The discrete adjoints~\eqref{eq:ODEadj_rec_d-o_euler} and~\eqref{eq:ODEadj_rec_d-o_mod_euler} from discretize-optimize converge to this same continuous adjoint, as one would expect. The key distinction is that discretize-optimize produces the \textit{exact} gradient of its respective discrete objective at any finite $\Delta t$, while optimize-discretize produces an approximation that becomes exact only in the limit.

As we will see next, unlike their deterministic counterpart, stochastic trajectories are not smooth (differentiable) in general, so what may appear as merely different discretizations of the same continuous model can yield different limits as the mesh gets finer.

\subsection{A Naive Attempt at Optimize-Discretize for It\^{o} SDEs}\label{subsec:naive_o-d_Ito}

Following our analysis of ODEs, let us now examine if the Euler-Maruyama discrete adjoint~\eqref{eq:ito_discrete_adj} converges to well-defined continuous adjoint process as $\Delta t \to 0$. 

The discrete adjoint~\eqref{eq:ito_discrete_adj} does not directly correspond to a continuous It\^{o} SDE in its current form. In a forward-time It\^{o} process, coefficients multiplying $dW(s)$ must depend only on information available at time $s$~\cite[Chap.~3.1]{oksendal2003} --- the defining left-endpoint convention. Here, the coupling of the Brownian increment $\Delta W_n$ with the future adjoint $p_{n+1}$ violates this. Interpreting the recursion backward in time does not help either as the Jacobian $\frac{\partial g(t_n,X_n)}{\partial X_n}$ evaluation at $(t_n, X_n)$ would then represent a right-endpoint evaluation coupled with $\Delta W_n$, again violating It\^{o}'s defining convention. 

It is tempting at this point to attempt to resolve this discrepancy by shifting all evaluations in the discrete adjoint from $(t_n, X_n)$ to $(t_{n+1}, X_{n+1})$, making the recursion resemble a backward-Euler step. However, this reasoning fails here and is misleading.
To see this, consider a backward-Euler accumulation of the adjoint-recursion \eqref{eq:ito_discrete_adj} along a sample path by shifting all evaluations from $(t_n,X_n)$ to $(t_{n+1}, X_{n+1})$
\begin{equation}\label{eq:ito_discrete_adj_o-d}
    \begin{cases}
        \widetilde{p}_n = \left( \mathbf{I} + \Delta t \frac{\partial f(t_{n+1},X_{n+1})}{\partial X_{n+1}}^\top + \Delta W_n^\top \frac{\partial g(t_{n+1},X_{n+1})}{\partial X_{n+1}}^\top\right) \  \widetilde{p}_{n+1} \\ 
        \widetilde{p}_N = \frac{\partial \Phi(X_N)}{\partial X_{N}}^\top.
    \end{cases}
\end{equation}

At first glance, this looks harmless, and even appears to have a valid backwards It\^{o}-type update. But this is exactly the problem. If we unroll the recursions into sums, we obtain
\begin{equation*}
    \widetilde{p}_0 = \frac{\partial \Phi(X_N)}{\partial X_{N}}^\top + \sum_{n=0}^{N-1} \Delta t\frac{\partial f(t_{n+1},X_{n+1})}{\partial X_{n+1}}^\top \widetilde{p}_{n+1} + \sum_{n=0}^{N-1} \Delta W_n^\top\frac{\partial g(t_{n+1},X_{n+1})}{\partial X_{n+1}}^\top \widetilde{p}_{n+1},  
\end{equation*}
for the backward-Euler in~\eqref{eq:ito_discrete_adj_o-d} and
\begin{equation*}
    p_0 = \frac{\partial \Phi(X_N)}{\partial X_{N}}^\top + \sum_{n=0}^{N-1} \Delta t\frac{\partial f(t_n,X_n)}{\partial X_n}^\top p_{n+1} + \sum_{n=0}^{N-1} \Delta W_n^\top\frac{\partial g(t_n,X_n)}{\partial X_n}^\top p_{n+1}, 
\end{equation*}
for the discrete adjoint in \eqref{eq:ito_discrete_adj}. Although the shift in evaluation points is generally benign for smooth trajectories, non-trivial errors can accumulate for non-smooth trajectories. Since stochastic paths are generically non-smooth, the two discretizations above will generally exhibit a non-zero discrepancy that will not vanish with the mesh refinement.

Hence, a naive attempt at optimize-discretize for It\^{o} SDEs yields a biased quantity and does not correctly recover the gradient we seek. Interestingly, when we compute the naive adjoint~\eqref{eq:ito_discrete_adj_o-d} for the Black--Scholes model, we find that it produces identical results to the discrete adjoint. 

This coincidence reflects a common phenomenon worth noting. Recall that the discrete adjoint recursion \eqref{eq:ito_discrete_adj} involves the Jacobian of the diffusion. For Black--Scholes with $g(S) = \sigma S$, we have $\frac{\partial g(S)}{\partial S} = \sigma$, which is constant. This means that evaluating the Jacobian at $(t_n, S_n)$ versus $(t_{n+1}, S_{n+1})$ yields the same value: $\frac{\partial g(S_n)}{\partial S_n} = \frac{\partial g(S_{n+1})}{\partial S_{n+1}} = \sigma$. Consequently, the correct and naive adjoints coincide.

This observation has practical significance. For many SDEs arising in applications, particularly those with state-independent or linear multiplicative noise, the diffusion Jacobian is state-independent. In such cases, practitioners can leverage existing ODE differentiation tools without concern. But this is not always the case; to reveal the importance of proper discretization, we turn to the Constant Elasticity of Variance (CEV) model~\cite{cox1976} of Cox: 
\begin{equation*}
    dS(t) = r S(t) \, dt + \sigma S(t)^\beta \, dW(t), \quad S(0) = S_0,
\end{equation*}
where $\beta \in \mathbb{R}$ controls how volatility scales with the stock price level. The CEV model is an extension of the Black--Scholes model that attempts to improve option pricing by relaxing the assumption of constant volatility. When $\beta = 1$, we recover Black--Scholes; when $\beta \neq 1$, the diffusion Jacobian becomes
\begin{equation*}
   \frac{\partial g(S)}{\partial S}  = \beta \sigma S^{\beta - 1},
\end{equation*}
which is now state-dependent. For our next experiment, we choose $\beta = 1.33$ to ensure substantial state-dependence. This particular value is arbitrary and chosen purely for illustration. We use the same option parameters as before to generate $10^5$ Monte Carlo paths with $\Delta t = 10^{-3}$ for computing the Greek $\Delta = \partial C / \partial S_0$ with our discrete and naive adjoints. We use identical Brownian increments for both methods to enable direct comparison.

\begin{figure}[htp]
  \centering
  \includegraphics[width=0.78\linewidth]{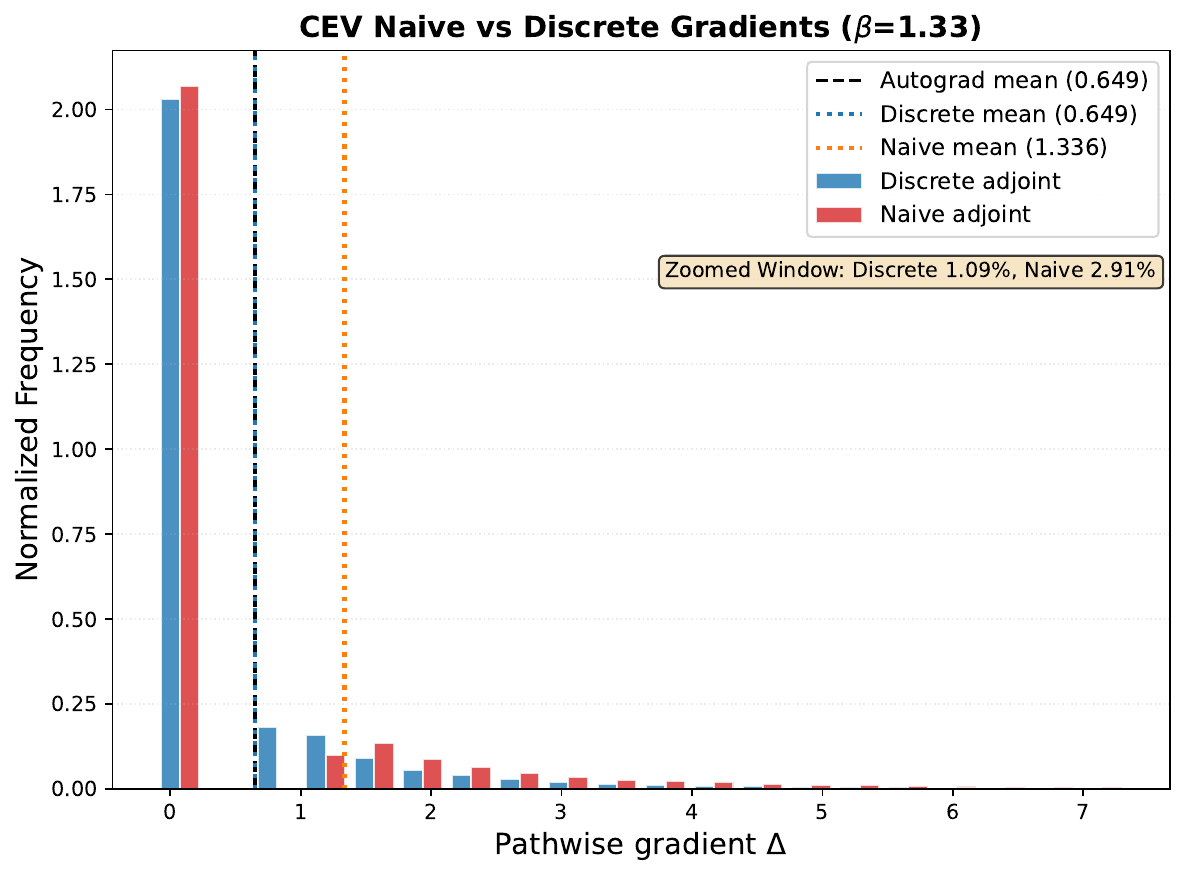}
  \caption{Histograms of pathwise gradients $\Delta = \partial C / \partial S_0$ for the CEV model with $\beta = 1.33$, comparing discrete-adjoint and naive-adjoint estimators. The x-axis is restricted to the pooled $[ 0.1\%, 98\% ]$ quantile range and the \enquote{Zoomed Window} note reports the per-method share of samples that fall beyond the right bound (Discrete $1.09\%$, Naive $2.91\%$). Computations using $10^5$ Monte Carlo paths with $\Delta t = 10^{-3}$.}
  \label{fig:naive_vs_discrete_adj}
\end{figure} 

\hyperref[fig:naive_vs_discrete_adj]{Figure~\ref*{fig:naive_vs_discrete_adj}} displays the histograms of the pathwise gradients. The discrete adjoint matches PyTorch's automatic differentiation as expected, both averaging $\Delta \approx 0.649$, while in contrast, the naive adjoint exhibits a systematic bias, with mean $\Delta \approx 1.336$ --- more than twice the true value. Moreover, the naive distribution has a heavy right tail, reflecting occasional explosive gradients.  This experiment illustrates clearly that when $\frac{\partial g(t,x)}{\partial x}$ is state-dependent, the discrete and naive approaches can diverge dramatically, with only the former producing unbiased gradient estimates. 

\subsection{Optimize-Discretize for Stratonovich SDEs} In~\Cref{subsec:naive_o-d_Ito}, we saw that a naive approach resembling the ODEs does not work for It\^{o} SDEs in general, but what about Stratonovich?

To this end, let us examine the backward recursion~\eqref{eq:heun_adjoint} more closely. Unrolling the recursion, we obtain
\begin{equation}\label{eq:unrolled_p_stratonovich}
   \begin{split}
       p_0 = & \ \frac{\partial \Phi(x_N)}{\partial X_N}^\top+\sum_{n=0}^{N-1} \frac{\Delta t}{2} \left[\frac{\partial f(t_n,X_n)}{\partial X_n}^\top + \frac{\partial f(t_{n+1},\widetilde{X}_{n+1})}{\partial \widetilde{X}_{n+1}}^\top\right]p_{n+1}\\
       &+ \sum_{n=0}^{N-1} \frac{\Delta W_n^\top}{2} \left[\frac{\partial g(t_n,X_n)}{\partial X_n}^\top + \frac{\partial g(t_{n+1},\widetilde{X}_{n+1})}{\partial \widetilde{X}_{n+1}}^\top\right] p_{n+1} + \sum_{n=0}^{N-1} H_n^\top p_{n+1},
   \end{split} 
\end{equation}
where $H_n$ is given by~\eqref{eq:higher_order_terms_heun}.

Focusing on the first two sums, observe that both the drift and diffusion Jacobians are evaluated as midpoint averages between $X_n$ and $\widetilde{X}_{n+1}$. Thus, when looking at~\eqref{eq:unrolled_p_stratonovich} in the time-reversed direction, the first two sums seem to mirror the structure of the Heun or trapezoidal schemes for approximating Stratonovich SDEs in a forward-time setting. Note that additional terms in $H_n$ arising from the predictor-corrector structure of the Heun scheme add some subtlety to the situation but do not alter this fundamental \enquote{reversible} character. 

Therefore, it is not unreasonable to expect that~\eqref{eq:unrolled_p_stratonovich} converges to a well-defined Stratonovich SDE as $\Delta t \to 0$. It turns out that this intuition is correct, and one can develop a rigorous theory of forward and backward Stratonovich SDEs based on the \textit{rough path theory}~\cite{friz2014, kidger2022neural, lyons1998}. The convergence delivered by this framework is \emph{pathwise}~\cite{lyons1998}, which is strictly stronger than weak convergence, and under sufficient integrability conditions it also implies strong convergence~\cite[Sec.~8.2--8.3]{kloeden1992}.

More specifically, one can show~\cite[Thm.~5.10]{kidger2022neural} that the gradient of the continuous objective function~\eqref{eq:objective_SDE_Ito} under Stratonovich dynamics~\eqref{eq:SDE_Stratonovich} is given by
\begin{equation}\label{eq:gradJ_cont_SDE_Stratonovich}
    \nabla_{x_0} J(x_0)=\mathbb{E} \left[p(0) \mid X_0=x_0\right],
\end{equation}
where the adjoint state $p$ satisfies the (suitably defined) adjoint SDE
\begin{equation}\label{eq:p_SDE_strat}
\begin{cases}
    dp(s)=-\frac{\partial f}{\partial x} (s,X(s))^\top p(s)ds-\sum_{j=1}^m\frac{\partial g_{:j}}{\partial x}(s,X(s))^\top p(s) \circ dW^j(s)\\
    p(T)=\frac{\partial \Phi (X(T))}{\partial x}^\top,
\end{cases}
\end{equation}
where $g_{:j}(s,x)=(g_{ij}(s,x))_{i=1}^d$ and $W(s)=(W^j(s))_{j=1}^m$.

Hence, one can use the optimize-discretize strategy to approximate the gradient~\eqref{eq:gradJ_cont_SDE_Stratonovich}, just as in the ODE setting, provided the numerical discretizations of~\eqref{eq:SDE_Stratonovich} and~\eqref{eq:p_SDE_strat} are chosen respectively so that they recover the Stratonovich SDE in both directions. Heun’s scheme is just one such choice, and this opens the door to higher-order numerical methods~\cite[Chap.~17]{higham2021} to achieve better convergence rates.

Parameters and running costs are handled exactly as in the discretize-optimize setting (see \Cref{subsec:params_running}): augmenting the state with $\theta$ (held constant), and with a running-cost coordinate when needed, leaves the dynamics a Stratonovich SDE, so the same continuous-adjoint construction applied to the augmented process yields the parameter gradient $\nabla_\theta J(x_0)$ alongside $\nabla_{x_0} J(x_0)$.

The conclusion is that the Stratonovich SDE renders a straightforward extension of the optimize-discretize approach to SDEs.

\begin{summarybox}[Summary: Continuous Adjoint for Stratonovich SDEs]
Assume that the state dynamics is given by the Stratonovich SDE
    \[
    \begin{cases}
        dX(s) = f(s,X(s)) ds + g(s,X(s)) \circ dW(s), \quad 0 \leq s \leq T, \\
        X(0) = x_0.
    \end{cases}
    \]
    Then the gradient of the objective function
    \[
    J(x_0)=\mathbb{E}\left[ \Phi(X(T)) \ | \ X(0) = x_0 \right]
    \]
    is given by
    \[
    \nabla_{x_0} J(x_0)=\mathbb{E} \left[p(0) \mid X(0)=x_0\right],
    \]
    where the adjoint state dynamics are given by
    \[
    \begin{cases}
    dp(s)=-\frac{\partial f}{\partial x} (s,X(s))^\top p(s)ds-\sum_{j=1}^m\frac{\partial g_{:j}}{\partial x}(s,X(s))^\top p(s) \circ dW^j(s)\\
    p(T)=\frac{\partial \Phi (X(T))}{\partial x}^\top.
\end{cases}
    \]
\end{summarybox} 

\subsection{Optimize-Discretize for It\^{o} SDEs. Go through Stratonovich}\label{subsec:ito_via_strat}

In situations where we do not have a choice for the underlying stochastic dynamics, and it has to be It\^{o}, we can still apply the It\^{o}-Stratonovich correction trick to take advantage of the optimize-discretize approach.

Indeed, a transition between It\^{o} and Stratonovich SDE can be performed via the following formula~\cite[Eq.~(6.1.3)]{oksendal2003}:
\begin{equation}\label{eq:ito<->str-multi-d}
    \begin{split}
        dX(s)=&f(s,X(s))ds+g(s,X(s)) d W(s) \\
        &\Updownarrow\\
        dX(s)=&\hat{f}(s,X(s))ds+g(s,X(s))  \circ d W(s)\\
        \text{where}\quad  \hat{f}_i(s,x)=&f_i(s,x)-\frac{1}{2}\sum_{j=1}^m \sum_{k=1}^d \frac{\partial g_{ij}(s,x)}{\partial x_k} g_{kj}(s,x),\quad 1\leq i \leq d.
    \end{split}
\end{equation}

Hence, given underlying It\^{o} dynamics~\eqref{eq:SDE_Ito}, we can rewrite it as an equivalent Stratonovich SDE via~\eqref{eq:ito<->str-multi-d} and proceed with optimize-discretize safely.

To validate this approach, we return to the CEV model and compare two methods for computing $\Delta = \partial C / \partial S_0$: the direct discrete adjoint via Euler-Maruyama on the It\^{o} SDE, and the corrected optimize-discretize approach where we first convert to Stratonovich form using~\eqref{eq:ito<->str-multi-d}, then apply Heun's scheme forward and discretize the backward continuous adjoint dynamics using Heun's scheme again. Using the same parameters as before $(S_0 = 100, K = 110, r = 0.05, \sigma = 0.2, T = 1, \beta = 1.33)$, we simulate $10^5$ Monte Carlo paths with step size $\Delta t = 10^{-3}$. For each approach, we use identical Brownian increments for both methods to enable direct comparison.

\begin{figure}[htp]
  \centering
  \includegraphics[width=0.8\linewidth]{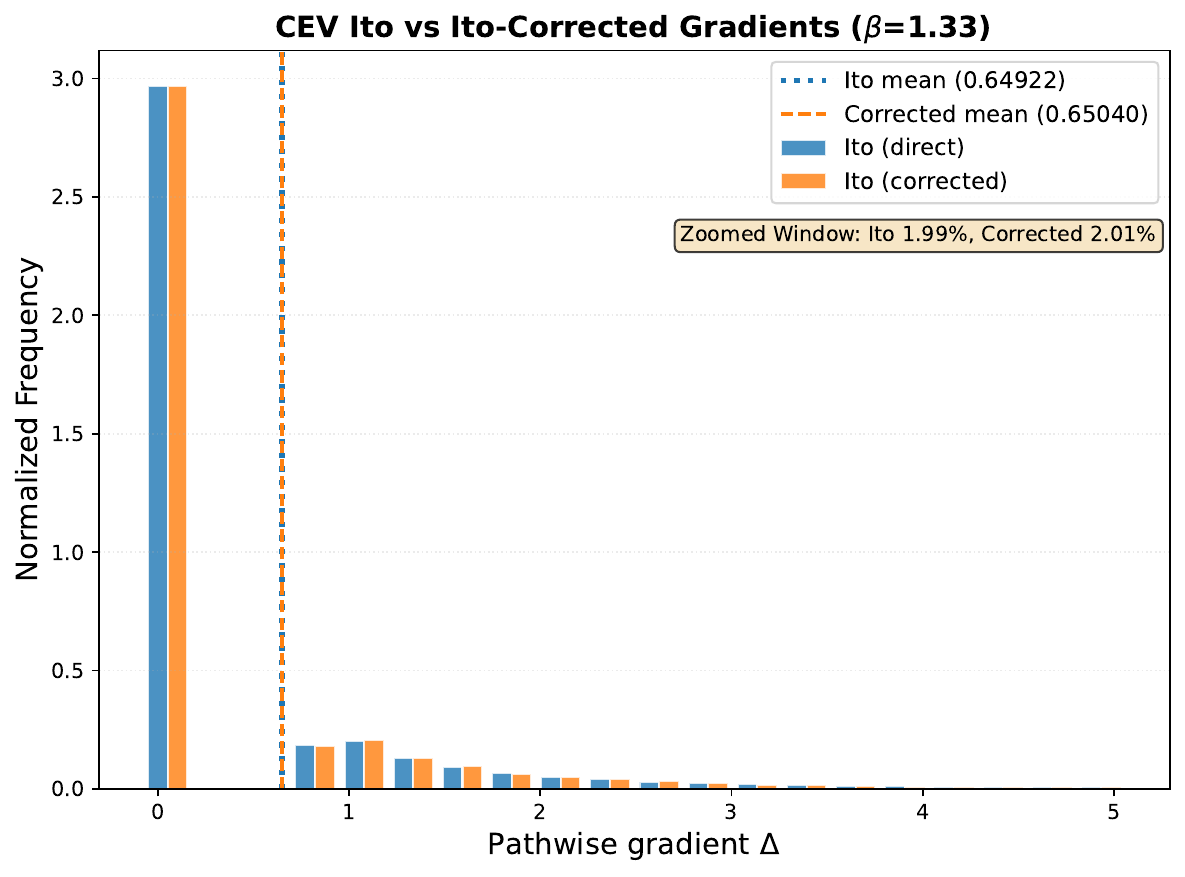}
  \caption{Histograms of pathwise gradients $\Delta = \partial C / \partial S_0$ for the CEV model with $\beta = 1.33$, comparing It\^{o} discrete adjoint (Euler-Maruyama) versus corrected Stratonovich adjoint (Heun). The x-axis is restricted to the pooled $[0.1\%, 98\%]$ quantile range, and the \enquote{Zoomed Window} note reports the per-method share of samples that fall beyond the right bound (It\^{o} $1.99\%$, Corrected $2.01\%$). Computations using $10^5$ Monte Carlo paths with $\Delta t = 10^{-3}$.}
  \label{fig:ito_strat_comparison}
\end{figure} 

\hyperref[fig:ito_strat_comparison]{Figure~\ref*{fig:ito_strat_comparison}} demonstrates that both approaches produce nearly identical gradient frequencies, validating the correction procedure between It\^{o} and Stratonovich frameworks. The slight discrepancy arises because the two approaches compute different adjoints: the direct It\^{o} via Euler-Maruyama discretization forward computes the exact discrete adjoint of its forward scheme, whereas the corrected It\^{o} via Heun discretization forward employs the optimize-discretize approach by discretizing the corresponding continuous Stratonovich adjoint equation with the Heun scheme, producing an approximation that converges to the exact gradient only as $\Delta t \to 0$. Nevertheless, this experiment demonstrates that both are unbiased estimators of the same continuous gradient in the limit, and both approximately agree.

\begin{summarybox}[Summary: Continuous Adjoint for It\^{o} SDEs]
    Given an It\^{o} state dynamics
    \[
    \begin{cases}
        dX(s) = f(s,X(s)) ds + g(s,X(s)) dW(s), \quad 0 \leq s \leq T, \\
        X(0) = x_0,
    \end{cases}
    \]
    pass to the corresponding Stratonovich dynamics
    \[
    \begin{cases}
        dX(s) = \hat{f}(s,X(s)) ds + g(s,X(s)) \circ dW(s), \quad 0 \leq s \leq T, \\
        X(0) = x_0.
    \end{cases}
    \]
    and compute the gradient of the objective function
    \[
    J(x_0)=\mathbb{E}[\Phi(X(T)) \ | \ X(0) = x_0]
    \]
    by the optimize-discretize approach as applied to the Stratonovich form.
\end{summarybox}

\subsection{The Lesson. Practical Implementation Considerations}

Having explored both discretize-optimize and optimize-discretize approaches, there are some practical insights to keep in mind when employing them. 

Discretize-optimize should use checkpointing when there is insufficient memory to store intermediate states for the backward pass. Rather than storing all intermediate states $X_n$, checkpointing stores states only at coarse intervals. During the backward pass, intermediate steps are regenerated by re-running the forward solver from the nearest checkpoint using the same Brownian increments. This is necessary because discretize-optimize computes the exact gradient of the discrete forward scheme, which requires the exact forward states at each time step.

Optimize-discretize can reduce the memory footprint dramatically by recomputing the states backward in time during the differentiation as long as the adjoint system is well-behaved.  After solving the forward SDE to obtain $X_N$, the backward adjoint equation is integrated backward in time together with the forward SDE as a coupled system using the reversed Brownian path. The forward SDE integrated backward with reversed increments approximately reconstructs the original forward trajectory due to the reversibility properties of Stratonovich integration~\cite[Rem.~5.12]{kidger2022neural}, \cite{friz2014, lyons1998}. This eliminates the need to store intermediate states $X_n$, trading computational cost (solving both equations backward) and approximation error (the reconstruction is approximate, not exact) for memory efficiency.

A key insight is that both approaches should utilize identical Brownian realizations for forward and backward passes. For discretize-optimize, this consistency is necessary to differentiate through the exact forward computation graph. For optimize-discretize, the same Brownian path enables the coupled backward integration described above.

Crucially, when using the same Brownian path for forward and backward passes, both methods compute unbiased per-sample gradients: for any single Brownian path, the pathwise gradient from that simulated state trajectory is an unbiased estimator of $\nabla_{x_0} J(t,x_0)$. This enables efficient optimization: simulate one forward path, compute its gradient, and average over many paths. Using different Brownian paths for forward and backward passes would break this property and introduce a bias, requiring Monte Carlo averaging over many Brownian paths during the backward pass for each forward trajectory, in addition to the standard averaging across forward trajectories. Matching paths preserves this essential per-sample efficiency.

\section{Practical Guidance and Outlook}

We now conclude by providing targeted guidance for readers from different communities and indicate where the techniques presented in this tutorial connect to broader frameworks and ongoing research.

\paragraph{Mathematics} For readers interested in a more abstract and unified treatment, controlled differential equations (CDEs) and rough path theory provide a framework that encompasses both ODEs and SDEs. In this setting, SDEs are viewed as special cases of CDEs driven by Brownian paths, and the pathwise approach to differentiation extends naturally to rough paths. This perspective also accommodates stochastic processes that are not classical SDEs, including fractional-Brownian and jump-driven dynamics, via the appropriate rough-path lift; the discretize-optimize lesson of Tale~1 carries over essentially unchanged whenever the forward scheme remains differentiable although the bias and variance of the resulting estimators depend on the scheme. For comprehensive mathematical treatments, see~\cite{ friz2014,kidger2022neural, lyons1998}. The techniques in this tutorial align directly with this rough path viewpoint, where backward integration along realized trajectories is the natural mode of computation.

\paragraph{Computational Mathematics} Differentiating through SDEs faces similar memory and computational challenges as ODEs, but with additional subtleties arising from the stochastic increments. For recent advances in scalable gradient computation for SDEs and regenerating Brownian paths see~\cite{jelin2024, kidger2021efficient}. For more of an introduction to memory-efficient techniques such as checkpointing for ODEs, see~\cite[Ch. 12]{griewank2008}. 

Additionally, the discretize-optimize approach naturally extends to higher-order numerical methods~\cite[Chap.~17]{higham2021} by differentiating through whatever forward scheme is chosen. On the other hand, the lessons of~\Cref{sec:tale_2} extend to higher-order methods as well: Stratonovich schemes (Milstein, Runge-Kutta) produce discrete adjoints that converge to well-defined continuous Stratonovich adjoints, enabling optimize-discretize strategies with improved accuracy; see~\cite{foster2023convergence, foster2024high} for additional higher-order schemes for Stratonovich SDEs. For It\^{o} schemes beyond Euler-Maruyama, the non-adapted nature of discrete adjoints persists, making discretize-optimize the preferred computational strategy.

Finally, we note that the adjoint techniques discussed in this tutorial extend to solving high-dimensional partial differential equations (PDEs). Deep BSDE methods~\cite{han2018solving, li2024neural} exploit the relationship between PDEs and SDEs to attempt to overcome the curse of dimensionality. In these works, neural networks approximate PDE solutions by reformulating the problem using a pair of SDEs that encode gradient information analogous to the adjoint processes derived here.

\paragraph{Optimal Control} Optimal control theory is one natural field where gradient (sensitivity) computations with respect to given parameters (controls) are required. Indeed, gradient-based optimization methods are among the most effective ways to compute optimal controls. Interestingly, many optimal control textbooks either omit discussing adjoints in the stochastic setting and employ the Hamilton-Jacobi methodology~\cite{fleming2006}, or they discuss \textit{adapted adjoints} which satisfy suitably augmented It\^{o} SDEs~\cite[Chap.~6]{pham2009}, \cite[Sec.~3.1]{yong1999}. More specifically, the adapted adjoint, $p^A(s)\in \mathbb{R}^d$, corresponding to the It\^{o} dynamics~\eqref{eq:SDE_Ito} and objective function
\[
J(s,x)=\mathbb{E}\left[\Phi(X(T))\mid X(s)=x \right]
\]
satisfies the SDE
\begin{equation} \label{eq:adapted_ito_adjoint}
    \begin{cases}
       
        dp^A_k(s) = \left[ - p^A(s) \cdot \partial_{x_k} f(s, X(s)) - \sum_{ij} q_{ij}(s) \partial_{x_k} g_{ij}(s, X(s))\right] ds  \\
        \qquad\qquad + \sum_j q_{kj}(s) dW_j(s),\quad 1\leq k \leq d, \\
        p^A(T) = \nabla \Phi(X(T)), 
    \end{cases}
\end{equation}
where $q(s) \in \mathbb{R}^{d \times m}$ is an auxiliary adapted process related to the Hessian of the objective function: $q(s) = \nabla_x^2 J(s, X(s)) g(s, X(s))$. The adapted adjoint $p^A$ is convenient due to the fact that, under suitable regularity conditions, one has that $p^A(s)=\nabla J(s,X(s))$, which yields optimal controls in a \textit{feedback form}~\cite[Sec.~6.4]{pham2009}, \cite[Chap.~3, Sec.~5]{yong1999}.

Unfortunately,~\eqref{eq:adapted_ito_adjoint} is not closed, and computing $p^A(s)$ requires simultaneously solving for $q(s)$, which depends on second derivatives unavailable from a single trajectory. This makes adapted adjoints considerably more expensive to approximate numerically, and they are often reserved for cases when~\eqref{eq:adapted_ito_adjoint} can be integrated analytically~\cite[Sec.~6.5]{pham2009}, \cite[Chap.~5]{yong1999}.

In contrast, the adjoints discussed in this paper are non-adapted, and do not require computing or estimating auxiliary processes (e.g.,~\eqref{eq:p_SDE_strat} is a closed system) and admit a seamless pathwise computation. Thus, the methods discussed in this paper can serve as versatile computational tools for stochastic optimal control problems. 

\paragraph{Finance} In quantitative finance~\cite{hull2018}, It\^{o} calculus is the dominant framework due to regulatory standards and the importance of the martingale property for no-arbitrage pricing. When continuous adjoints are required, practitioners should employ a double correction procedure. Starting with an It\^{o} SDE, one must first convert to Stratonovich form via the drift correction, derive the corresponding Stratonovich discrete adjoint equation, and then convert the adjoint back to It\^{o} form for interpretation within the standard financial mathematics framework. For practitioners who use Monte Carlo methods~\cite[Chap.~7]{glasserman2003} to compute Greeks but do not require continuous adjoints, the discretize-optimize approach with a proper numerical discretization of It\^{o} offers a more direct alternative that avoids these multiple conversions. For non-smooth payoffs such as digital options, the pathwise approach used throughout this primer is biased; the standard remedy is either to smooth the payoff or to switch to a likelihood-ratio or Malliavin-calculus estimator that differentiates the law of $S_T$ rather than the path~\cite{fournie1999}.

\paragraph{Machine Learning} Differentiating through stochastic differential equations has become central to modern machine learning. Neural ordinary differential equations~\cite{chen2018} and their stochastic extensions~\cite{kidger2022neural} parameterize dynamics with neural networks and require differentiating through the dynamics to train. Additionally, CDEs provide continuous-time analogues of recurrent neural networks~\cite{kidger2020neural}, enabling architectures for irregular time series and sequential data.  For neural SDEs and neural CDEs, Stratonovich calculus should be preferred precisely because there are no martingale constraints or measure-theoretic requirements --- the focus is purely on pathwise dynamics and their gradients. The techniques in this tutorial directly apply to training such models.

As stochastic modeling continues to expand across scientific computing, finance, and machine learning, the ability to correctly compute gradients of stochastic systems remains essential. This tutorial has sought to equip readers with the essential techniques for computing sensitivities of stochastic differential equations. For practitioners seeking to differentiate their stochastic system, the simplest approach is to write the forward SDE solver in a framework that enables automatic differentiation, such as PyTorch or JAX, and apply standard backpropagation.  Armed with the principles and the intuition developed through our two-tale narrative, we hope that readers can now approach gradient computation for SDEs with confidence.

\section*{Acknowledgments} 

We are very grateful to James Foster, Patrick Kidger, and Richard Lehoucq for insightful discussions, helpful pointers to the literature, and careful reading and feedback on an earlier version of this manuscript.

We used large language models (LLM) to polish the manuscript for clarity. Claude Code assisted with implementing and testing the numerical experiments to improve quality and ensure reproducibility. We take full responsibility for all content.

\end{document}